\documentclass[12pt,a4paper]{amsart} 
\usepackage{enumerate}

\ifx\pdftexversion\undefined
\usepackage[dvips]{graphicx,color}
\else
\usepackage[pdftex]{graphicx,color}
\fi

\let\ifanglais\iftrue

\def\R{{\mathbb R}}
\def\N{{\mathbb N}}
\def\pP{{\mathbb P}}
\def\vol{{\rm Vol}}
\def\aire{{\rm Area}}

\newcommand{\cC}{\mathcal{C}}
\newcommand{\cG}{\mathcal{G}}

\newcommand{\cS}{\mathcal{S}}
\newcommand{\cD}{\mathcal{D}}

\newcommand{\imp}{\Longrightarrow} 

\newcommand{\disp}{\displaystyle}
\newcommand{\seq}[2]{\left( #1_{#2} \right)_{\! #2 \in \N}} 
\newcommand{\bC}{\partial \cC}
\newcommand{\bD}{\partial \cD}

\newcommand{\bS}{\partial \cS}

\renewcommand{\bar}{\overline}
\newcommand{\ed}[1]{\textrm{d} #1}

\newtheoremstyle{mesthm}
  {10pt plus 1pt minus 1pt}
  {9pt plus 1pt minus 6pt}
  {\slshape}
  {0.5cm}
  {\bfseries}
  {.}
  {1ex}
  {}
\newtheoremstyle{mesdefi}
  {6pt plus 1pt minus 1pt}
  {6pt plus 1pt minus 1pt}
  {}
  {0.5cm}
  {\bfseries}
  {.}
  {1ex}
  {}%

\theoremstyle{mesthm}
\newtheorem{lema}{\ifanglais{\large L}emma\else{\large L}emme\fi}
\newtheorem{theo}[lema]{\ifanglais{\large T}heorem\else {\large
    T}h\'eor\`eme\fi}

\newtheorem{prop}[lema]{{\large P}roposition}  
\newtheorem{cor}[lema]{{\large C}orollary}

\newtheorem{rmq}[lema]{\ifanglais{\large R}emark\else{\large
    R}emarque\fi}

\theoremstyle{mesdefi}
\newtheorem{defi}[lema]{\ifanglais{\large D}efinition\else{\large
    D}\'efinition\fi}

\newcounter{step}
\renewcommand{\thestep}{\arabic{step}}
\newcommand{\step}{\refstepcounter{step}
\vskip 9pt plus 1pt minus 1pt
\par \indent 
{\bf \hspace{0.2em} Step \thestep:}
}
\newcommand{\istep}{\setcounter{step}{0}}

\newenvironment{enumebis}{\begin{itemize}\item[]\begin{enumerate}}{\end{enumerate}\end{itemize}}

\title[Ideal triangles and $\delta$-hyperbolicity in Hilbert Geometry]{Area of ideal triangles and Gromov hyperbolicity
in Hilbert Geometry}
\author{B.~Colbois, C.~Vernicos  and  P.~Verovic}

\begin{document}
\maketitle

\section*{Introduction and statements}

The aim of this paper is to show, in the context of Hilbert geometry, the equivalence between the
existence of an upper bound on the area of ideal triangles and the Gromov-hyperbolicity.

\begin{figure}[htpb] 
  \centering
  \includegraphics[scale=.7]{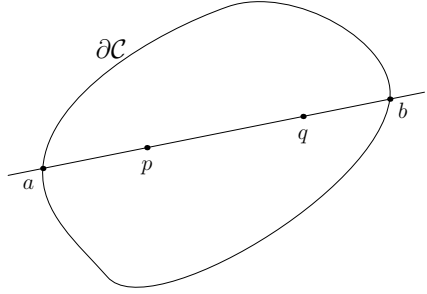}
  \caption{The Hilbert distance \label{dintro}}
\end{figure}

Let us recall that a Hilbert geometry
$(\mathcal{C},d_\mathcal{C})$ is a non empty bounded open convex set $\mathcal{C}$
on $\R^n$ (that we shall call \textit{convex domain}) with
the Hilbert distance 
$d_\mathcal{C}$ defined as follows : for any distinct points $p$ and $q$ in $\mathcal{C}$,
the line passing through $p$ and $q$ meets the boundary $\partial \mathcal{C}$ of $\mathcal{C}$
at two points $a$ and $b$, such that one walking on the line goes consecutively by $a$, $p$, $q$
$b$ (figure~\ref{dintro}). Then we define
$$
d_{\mathcal C}(p,q) = \frac{1}{2} \ln [a,p,q,b],
$$
where $[a,p,q,b]$ is the cross-product of $(a,p,q,b)$, i.e., 
$$
[a,p,q,b] = \frac{\| q-a \|}{\| p-a \|} \times \frac{\| p-b \|}{\| q-b\|} > 1,
$$
with $\| \cdot \|$ the canonical euclidean norm in
$\mathbb R^n$.

\smallskip
Note that the invariance of the cross-product by a projective map implies the invariance 
of $d_{\mathcal C}$ by such a map.

\smallskip
These geometries are naturally endowed with
a  $C^0$ Finsler metric $F_\mathcal{C}$ as follows: 
if $p \in \mathcal C$ and $v \in T_{p}\mathcal C =\R^n$
with $v \neq 0$, the straight line passing by $p$ and directed by 
$v$ meets $\partial \mathcal C$ at two points $p_{\mathcal C}^{+}$ and
$p_{\mathcal C}^{-}$~; we then define
$$
F_{\mathcal C}(p,v) = \frac{1}{2} \| v \| \biggl(\frac{1}{\| p -
  p_{\mathcal C}^{-} \|} + \frac{1}{\| p - p_{\mathcal C}^{+}
  \|}\biggr) \quad \textrm{et} \quad F_{\mathcal C}(p , 0) = 0.
$$ 

\begin{figure}[htpb]
  \centering 
    \includegraphics[scale=.7]{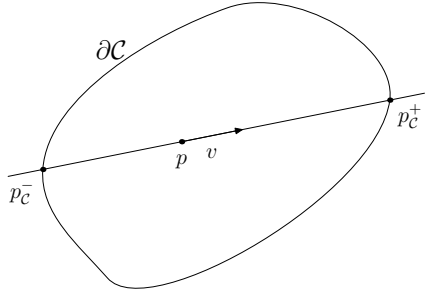}
          \caption{The Finsler structure \label{dintro2}}
          \end{figure}

\smallskip
The Hilbert distance $d_\mathcal{C}$ is the length distance associated to 
$F_{\mathcal C}$.

\smallskip
Thanks to that Finsler metric, we can built a Borel measure
$\mu_{\mathcal C}$ on $\mathcal C$ (which is actually
the Hausdorff measure of the metric space $(\mathcal C,
d_{\mathcal C})$, see \cite{bbi}, exemple~5.5.13 ) as follows.

\smallskip
To any $p \in \mathcal C$, let $B_{\mathcal C}(p) = \{v \in \R^n ~|~ F_{\mathcal{C}}(p,v) < 1 \}$
be the open unit ball in
$T_{p}\mathcal{C} = \R^n$ of the norm $F_{\mathcal{C}}(p,\cdot)$ and 
$\omega_{n}$ the euclidean volume of the open unit ball of the standard euclidean space
$\R^n$.
Consider the (density) function $h_{\mathcal C}\colon  \mathcal C \longrightarrow \R$ given by $h_{\mathcal C}(p)
= \omega_{n}/\vol\bigl(B_{\mathcal C}(p)\bigr),$ where $\vol$ is the canonical Lebesgue measure
of $\R^n$. We define $\mu_{\mathcal C}$, which we shall call
the \textit{Hilbert Measure} on $\mathcal{C}$,
by
$$
\mu_{\mathcal C}(A) = \int_{A} h_{\mathcal C}(p) \ed{\vol(p)}
$$
for any Borel set $A$ of $\mathcal C$.

\medskip
A fundamental result of Y. Benoist \cite{benoist} 
gives an extrinsic characterization of Gromov-hyperbolic
Hilbert geometries, that is sufficient and necessary conditions on the boundary $\partial \mathcal C$ of
a convex domain $\mathcal C$ to insure that the associate Hilbert geometry $(\mathcal C,d_{\mathcal C})$ is
Gromov-hyperbolic.

\medskip
The goal of this paper is to give an intrinsic condition equivalent to the Gromov-hyperbolicity
in terms of the area of the ideal triangles of $(\mathcal C,d_{\mathcal C})$.

\medskip
We define an ideal triangle $T \subset \mathcal C$ as the affine convex hull of three points $a, b, c$
of $\partial \mathcal C$ not on a line, and such that $T\cap \partial \mathcal C = 
{a} \cup {b} \cup{c}$. (Note that the affine convex hull coincide with the geodesic
convex hull when the space is uniquely geodesic, which is the case of Gromov-hyperbolic
Hilbert geometry). The area of a triangle $T$ (ideal or not) of $(\mathcal C,d_{\mathcal C})$,
denoted by $\aire_{\mathcal{C}}(T)$,
is its area for the Hilbert measure of $(\mathcal C \cap P , d_{\mathcal C \cap P})$,
where $P$ is the unique plane in $\R^n$ containing the triangle (in dimension $2$, as 
$\mathcal C \cap P=\mathcal C$, we will also denote it by $\mu_{\mathcal C}(T)$).  

\bigskip
In this paper, we prove

\begin{theo}
Let $\delta >0$. There exist a constant $C(\delta)>0$ with the following property: 
the Hilbert geometry $(\mathcal C,d_{\mathcal C})$ is $\delta$-hyperbolic if and
only if the area of any ideal triangle $T\subset \mathcal C$
is bounded above by $C(\delta)$.
\end{theo}

\bigskip
To show that the bound on the area of ideal triangles implies the $\delta$-hyperbolicity
is quite straightforward and its proof is in the first part of the paper (Theorem \ref{volIdelta}). The converse is
much more delicate: we show it on the second part of the paper (Theorem \ref{thmhyp-bounded}). The main ingredient
of the proof is a co-compacity Lemma (Theorem \ref{lemmacocompactness}, whose idea goes back in some sense to Benzecri \cite{benzecri})
and the results of Benoist's paper \cite{benoist}. To make the proof readable, 
we let some technical lemma in an Appendix at the end of the paper, in particular the Lemma \ref{lemHolder} 
deduced from \cite{benoist}, which implies an $\alpha$-H\"older regularity of the boundary of a convex 
domain whose Hilbert geometry is $\delta$-hyperbolic, with $\alpha$ depending only on $\delta$, and 
Lemma \ref{lemfinite-volume}, where we show that the $\alpha$-H\"older regularity implies the
finiteness of the area of ideal triangles.

\smallskip
Note that the results of this Appendix are used many times in the proof 
of Theorem \ref{thmhyp-bounded}.

\medskip
In the sequel we will 
switch between affine geometry (where our results are stated) and projective geometry (where Benoist's results are stated).
We will use the following two classical facts (see \cite{ps} section 1.3 page 8--11)
\begin{enumerate}
\item Any affine space can be embedded into a projective space
(by "adding an hyperplane at infinity"). Furthermore any one-to-one affine map extends
to a homography keeping the "hyperplane at infinty" globally invariant.
\item The complement of a projective hyperplane in a projective space is an affine
space. Furthermore all homographies keeping this hyperplane globally invariant
are naturally identified with an affine map on the complement.
\end{enumerate}

\section{Bounded area implies $\delta$-hyperbolicity}
In this part we prove
\begin{theo}\label{volIdelta} Let $M>0$. There exists $\delta=\delta(M)>0$ with the following property:
Let $(\mathcal{C},d_\mathcal{C})$ be a convex domain with its induced Hilbert distance.
If any ideal triangle in $(\mathcal{C},d_\mathcal{C})$ has its area less
than $M$ then $(\mathcal{C},d_\mathcal{C})$ is $\delta$-hyperbolic.
\end{theo}

This theorem is a straightforward consequence of the following proposition:

\begin{prop}\label{propimp} 
There exist a constant $C>0$ with the following property: for any $\delta$,
  if  $(\mathcal{C},d_\mathcal{C})$ is not $\delta$-hyperbolic, then there exists
an ideal triangle $T \subset \mathcal C$, whose area satisfies $\mu_{\mathcal C}(T) \geq C\cdot\delta$.
\end{prop}

Indeed, if the assumption of Theorem \ref{volIdelta} are satisfied, then $\mathcal{C}$
has to be $\delta$-hyperbolic for any $\delta>M/C$, otherwise we would get a contradiction
with the Proposition \ref{propimp}.

\smallskip
Now let us prove Proposition \ref{propimp}.
We already know that if $\partial \mathcal C$ is not strictly convex, then there is an ideal
triangle of arbitrarily large area (\cite{cvv} Corollaire 6.1 page 210).
Hence we can assume that $\partial \mathcal C$ is strictly convex, which 
implies that all the geodesics of $(\mathcal{C},d_\mathcal{C})$ are straight segments
(see \cite{dlharpe} proposition 2 page 99). 

Each
triangle $T \subset \mathcal C$ determines a plane section of $\mathcal C$, 
and is contained in an ideal triangle of this plane section. 
So, it suffices to exibit
a triangle (not necessarily ideal) such that  $\mu_{\mathcal C}(T) \geq C\cdot\delta$.

This is done thanks to the two following lemma.
\begin{lema}\label{lematd}
  If $(\mathcal{C},d_\mathcal{C})$ is not $\delta$-hyperbolic, there is a plane $P$
and a triangle $T$ in $P\cap\mathcal{C}$ such that a point in the triangle is
at a distance greater than $\delta/4$ from its sides.
\end{lema}

\begin{proof}[Proof of lemma \ref{lematd}]
If $(\mathcal{C},d_\mathcal{C})$ is not $\delta$-hyperbolic, there
exists a triangle $T \in \mathcal C$ of vertices $a,b,c$, a point $p \in \partial T$, say between $a$ and $b$,
 such that the distance from $p$ to the
two opposite sides of $\partial T$ is greater than $\delta$. 
The end of the proof takes place in the plane determined by the triangle $T$.

\smallskip
Let $R=\delta/2$. Consider a circle
$S$ of center $p$ and radius $R$. Let $p_1,p_2= S\cap \partial T$. We have $d_\mathcal{C}(p_1,p_2)=2R$.

\smallskip
If $q\in S$, then $d_\mathcal{C}(p_1,q)+d(q,p_2)\geq 2R$ by the triangle inequality. 
By continuity, we can choose $q\in S \cap T$, 
with $d_\mathcal{C}(q,p_1)\geq R;\ d_\mathcal{C}(q,p_2)\geq R$. 
From this fact, and by the classical triangular 
inequality, we deduce $d_\mathcal{C}(q,\partial T)\geq R/2$: to see it, 
let $p_3$ be the middle of the segment $pp_1$. We have
$d_\mathcal{C}(p,p_3)=d_\mathcal{C}(p_3,p_1)=R/2$. 

\begin{itemize}
\item If $q'\in pp_3$, $d_\mathcal{C}(q,q') \geq d_\mathcal{C}(p,q)-d_\mathcal{C}(p,p_3) \geq R/2$.
\item  If
$q'\in p_3p_1$, then $d_\mathcal{C}(q,q')\geq d_\mathcal{C}(q,p_1)-d_\mathcal{C}(q'q_1) \geq R/2$ and this show also that if
$d_\mathcal{C}(q',p_1)\leq R/2$ then $d_\mathcal{C}(q,q') \geq R/2$.
\item If $q'$ is such that $d_\mathcal{C}(q',p) \geq 3R/2$, then
$d_\mathcal{C}(q,q')\geq R/2$.  
\end{itemize}

This allow to conclude for the half line issue from $p$ 
through $p_1$ and we can do the same for the other half line.
\end{proof}

\begin{lema}\label{arealowerbd}
  There exists a constant $C_n$ such that any ball of radius $R>2$ in any  Hilbert geometry
of dimension $n$ has a volume greater or equal to $C_n\cdot R$
\end{lema}

\begin{proof}
  Let $B$ a ball centered at $q$ of radius $R$.
Consider a geodesic segment starting at $q$: it has length $R$ and lies inside $B$. 
We can cover it by $N=\text{integer part of } R$, parwise disjoint balls of radius $1$
contained in $B$, with $N\to \infty$ with $\delta$. But we know (Theorem 12, \cite{cv1}) 
that the volume of a radius $1$
ball is uniformly bounded below for all the Hilbert geometries by a constant $c(n)$.
Hence the volume of the ball of radius $R\geq2$ is greater than $(R-1)\cdot c(n)\geq R\cdot c(n)/2$.
\end{proof}

Hence, if $(\mathcal{C},d_\mathcal{C})$ is not $\delta$ hyperbolic thanks to lemma \ref{lematd} we would
find a triangle $T$  containing a two-dimensional ball of radius $\delta/4$, hence its area would
be greater than $\delta/4\cdot C_2$ thanks to lemma \ref{arealowerbd}. Which ends the proof of proposition
\ref{propimp}. 

\medskip

A consequence of Theorem \ref{volIdelta}, already proved with different approaches by 
A.~Karlsson and G.A.~Noskov \cite{kn}, Y.~Benoist \cite{benoist} and
 B.~Colbois and P.~Verovic \cite{cv},  is the following:
\begin{cor}
If the boundary of $\mathcal{C}$ is $C^2$ with strictly positive curvature, then
 $(\mathcal{C},d_\mathcal{C})$ is Gromov hyperbolic.
\end{cor}
This is a consequence of Theorem 4 in \cite{cvv} which says that if the boundary is
$C^2$ with strictly positive curvature, then the assumptions of Theorem \ref{volIdelta} are satisfied.

\section{From $\delta$-hyperbolicity to bounded area}

The aim of this section is to prove the following

\begin{theo} \label{thmhyp-bounded}
Let $\delta >0$. Then, there exists $V=V(\delta)>0$ with the following property:   
Let $\mathcal{C}$ be a convex domain such that $(\mathcal{C},d_{\mathcal{C}})$ 
is $\delta$-hyperbolic. 
   Then, for any ideal triangle T of $\mathcal{C}$, we have
$\aire_{\mathcal{C}}(T)\leq V$.
\end{theo}

Thought the ideas to prove Theorem \ref{thmhyp-bounded} are quite simple, the proof itself is somewhat
technical. The bound on the area of ideal triangle depends only on the $\delta$ of
the Gromov hyperbolicity. Therefore it suffices to prove Theorem \ref{thmhyp-bounded} in the
two dimensional case. Thus, from this point on, everything will be done
in the two dimensional case.

\subsection{Co-compactness of triangle-pointed convex}
Le us begin with some notations.

\smallskip
Let $G_n:={\rm PGL}(\R^{n})$, $\pP^n:=\pP(\R^{n+1})$ the projective space of $\R^{n+1}$.
A properly convex subset $\Omega$ of $\pP^n$ is a convex set such that there
is a projective hyperplane who doesn't meet its closure.
Denote by $X_n$ is the set of properly convex open sets.
Let $X_n^\delta$ be the set of $\delta$-hyperbolic properly convex open sets in $\mathbb{P}^n$ 

\smallskip
In $X_n$ we will consider the topology induced by the 
Hausdorff distance between sets, denoted by $d$.

\smallskip
We will say that a convex domain $\mathcal{C}$ is \textit{triangle-pointed} if one
fixes an ideal triangle in $\mathcal{C}$. Let 
$$T_2^\delta=
\{(\omega,x,y,z)\in X_2^\delta\times\pP^2\times\pP^2\times\pP^2\mid x,y,z \in \partial \omega, x\neq y, y\neq z, z\neq x\} $$ 
be the set of triangle-pointed
convex sets $\mathcal{C}$ with $\mathcal{C}\in X_2^\delta$. 

One of the main steps of our proof will rely upon the following cocompactness
result.

\begin{theo}\label{lemmacocompactness}
  $G_2$ acts cocompactly on $T_2^\delta$, i.e., for any sequence $(\mathcal{\omega}_n,\Delta_n)_{n\in\N}$ in $T_2^\delta$,
there is a sequence $(g_n)_{n\in\N}$ in $G_2$ and a subsequence of $(g_n\omega_n,g_n\Delta_n)_{n\in\N}$ 
that converges to $(\omega,\Delta)\in T_2^\delta$. 
\end{theo}

Actually, Theorem \ref{lemmacocompactness} is a corollary of the following more precise statement

\begin{prop}\label{propcocompactness}
  Let $(\omega_n,T_n)_{n\in\N}$ be a sequence in $T_2^\delta$, then 
  \begin{enumerate}
  \item There is a sequence $(g_n)_{n \in \N}$ in $ G_2$ and a number $0<e\leq1/2$
such that $g_nT_n=\Delta\subset \R^2$ the triangle whose coordinates are the points $(1,0)$, $(0,1)$, $(1,1)$,
and $g_n\omega_n\subset \R^+\times\R^+$ is tangent at $(1,0)$ to the $x$-axe, at $(0,1)$ to the $y$-axe
and at $(1,1)$ to the line passing through the points $(1/\alpha_n,0)$ and $(0,1/(1-\alpha_n)$ for
some $0<e\leq\alpha_n\leq 1/2$;
\item From the previous sequence we can extract a subsequence converging to some
 $(\omega,\Delta)\in T_2^\delta$.
  \end{enumerate}

\end{prop}
\begin{proof} \textbf{Step 1: A first transformation}

\medskip
According to \cite{cvv} (Proof of Th\'eor\`eme~3, p.~215 and Lemme~9, p.~216), for each $n \in \N$, 
there exist a number $\alpha_{n} \in (0 , 1/2]$ and an affine transformation $A_{n}$ of $\R^2$ such that:

\begin{enumerate}[1)]  
   \item The bounded open convex domain $\Omega_{n} := A_{n}(\omega_n)$ is contained in the triangle 
$\mathcal{T} \subset \R^2$
   whose vertices are the points $(0 , 0)$, $(1 , 0)$ and $(0 , 1)$.
   
   \item The points $(\alpha_{n} , 0)$, $(0 , 1 - \alpha_{n})$ and $(\alpha_{n} , 1 - \alpha_{n})$ are in $\partial \Omega_{n}$
   and the ideal triangle $\Delta_{n}$ they define in $(\Omega_{n} , d_{\Omega_{n}})$ is equal to $A_{n}(T_{n})$.
    
   \item The $x$-axis, the $y$-axis and the line passing through $(1 , 0)$ and $(0 , 1)$ 
   are tangent to $\partial \Omega_{n}$ at the points $(\alpha_{n} , 0)$, $(0 , 1 - \alpha_{n})$ and $(\alpha_{n} , 1 - \alpha_{n})$ respectively.
\end{enumerate}

Remark that we may have to take out different projective lines
to see the proper convex sets $\omega_n$ as convex sets in an affine space. But
up to some homography we can suppose that we always took the same. The geometries
involved will not be changed.

\begin{figure}[hbtp]
  \centering
  \includegraphics{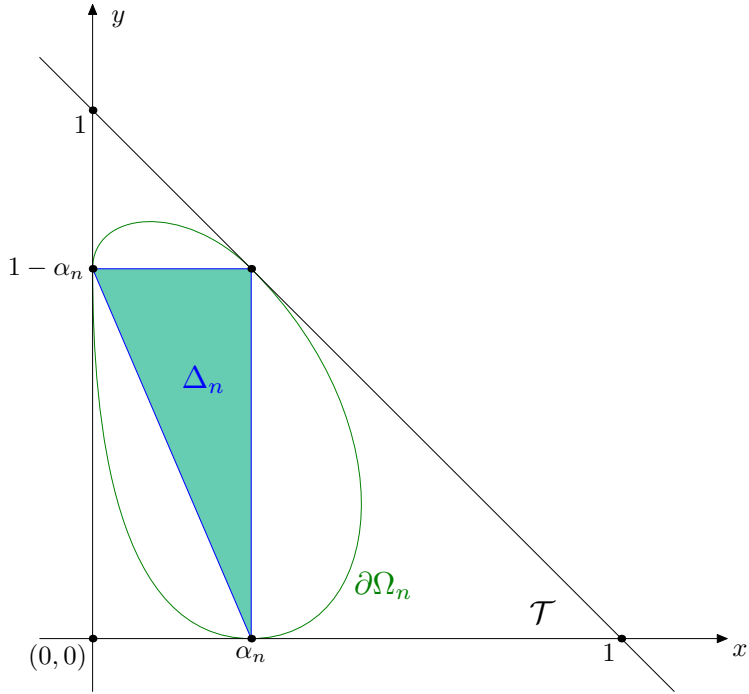}
  \caption{The $\Omega_n$ are convex sets included in a fixed triangle}
  \label{fig2}
\end{figure}

\vskip 9pt plus 1pt minus 1pt

\smallskip
\noindent
\textbf{Step 2: Proof of the first part of (1)} 

\medskip
In this part, we show the first part of point (1). The second point of (1),
that is to see that the set of $\{\alpha_n\}$ is uniformly bounded below by $e >0$, will be done at the step 4.

\medskip
For each $n \in \N$, if we consider the linear transformation $L_{n}$ of $\R^2$ defined by 
$$
L_{n}(1 , 0) = (1 / \alpha_{n} , 0)
\quad \mbox{and} \quad 
L_{n}(0 , 1) = (0 , 1 / (1 - \alpha_{n})),
$$
we have: 

\begin{enumerate}[1)]  
   \item The bounded open convex domain $\cC_{n} := L_{n}(\Omega_{n})$ is contained 
   in the triangle $\mathcal{T}_{n} \subset \R^2$
   whose vertices are the points $(0 , 0)$, $(1 / \alpha_{n} , 0)$ and $(0 , 1 / (1 - \alpha_{n}))$.
   
   \item The points $(1 , 0)$, $(0 , 1)$ and $(1 , 1)$ are in $\bC_{n}$
   and the ideal triangle $\Delta$ they define in $(\cC_{n} , d_{\cC_{n}})$ is equal to $L_{n}(\Delta_{n})$.
    
   \item The $x$-axis, the $y$-axis and the line passing through $(1 / \alpha_{n} , 0)$ and $(0 , 1 / (1 - \alpha_{n}))$ 
   are tangent to $\bC_{n}$ at the points $(1 , 0)$, $(0 , 1)$ and $(1 , 1)$ respectively.
\end{enumerate}

\begin{figure}[hbtp]
  \centering
  \includegraphics[scale=.7]{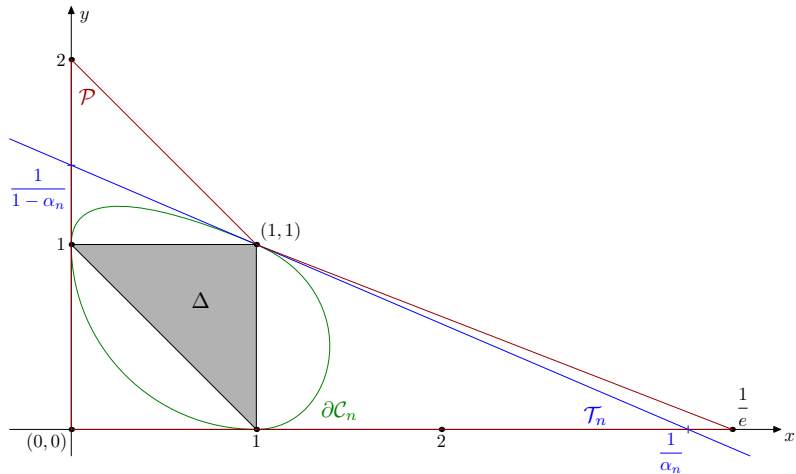}
  \caption{The $\mathcal{C}_n$ are convex sets with a fixed ideal triangle}
  \label{fig3}
\end{figure}

For each $n \in \N$, the affine transformation $L_{n} \circ A_{n}$ of $\R^2$ induces an isometry 
from $(\omega_n , d_{\omega_n})$ onto the metric space $(\cC_{n} , d_{\cC_{n}})$, 
Hence we have that $(\cC_{n} , d_{\cC_{n}})$ is $\delta$-hyperbolic for all $n \in \N$.

\smallskip
\noindent
\textbf{Step 3: Convergence of a subsequence} 

\medskip
All the convex domains $\cC_{n} \subset \R^2$ contain the fixed triangle $\Delta$ and are by construction
contained in the convex subset $\mathcal B =\{(x,y)\in \R^2: x \geq 0; 0\leq y\leq 2\}.$ The convex
$\mathcal B$ correspond to a properly convex set of the projective plane, because it does not contain the line $\{x=-1\}$.

From Lemma 2.2 page 189 in \cite{benoist}, 
the set of all the bounded open convex domains in the projective plane $\pP^2$ 
contained in $\mathcal{B}$ and containing the image of $\Delta$ 
is compact for the Hausdorff distance $d$. 
Thus there exist a proper convex domain $\Omega$ in $\pP^2$
such that $\Omega \subset \mathcal{B}$
and a subsequence of $\seq{\mathcal{C}}{n}$, still denoted by $\seq{\mathcal{C}}{n}$, 
such that $d(\mathcal{C}_n , \Omega) \to 0$ as $n \to +\infty$.

\smallskip
Point a) of Proposition~2.10, page 12, in Benoist \cite{benoist} then implies 
that $\Omega$ is $\delta$-hyperbolic and strictly convex. 

\smallskip
Note that since the points $(1 , 0)$, $(0 , 1)$ and $(1 , 1)$ are in $\bC_{n}$ for all $n \in \N$,
they also are in $\partial \Omega$.

\smallskip
\noindent
\textbf{Step 4: The bound on the $\alpha_n$} 

\medskip
\textbf{By contradiction:}  Suppose $\inf{\{ \alpha_{n} : n \in \N \}} = 0$.

\smallskip
By considering a subsequence, we can assume that 
$$ 
\lim_{n \to +\infty} \alpha_{n} = 0 \text{.}
$$

Then we have that for any $\mathcal{C}_n$, a part of its boundary is
in the triangle $(0,1)$, $(1,1)$ and $(0,1/(1-\alpha_n))$. When $n\to +\infty$, the last
point converges towards $(0,1)$, i.e. the triangle collapses on the segment
defined by $(0,1)$ and $(1,1)$. Hence, this segment is on $\partial \Omega$, which contradicts
the strict convexity of step $3$.

This implies that there exists a constant $e > 0$ such that $\alpha_{n} \in [e , 1/2]$ for all $n \in \N$, and
that $\Omega$ is bounded in $\R^2$.

\end{proof}


\begin{prop}\label{finitearea}
  Let $\mathcal{C}$ be a bounded open convex domain in $\R^2$ such that $\partial \mathcal{C}$ is 
$\alpha$-H\"older for some $\alpha>1$. 
Then for any ideal 
triangle $T$ in $(\cC,d_{\mathcal{C}})$, $\mu_{\mathcal{C}}(T)$ is finite.
\end{prop}

\begin{proof}
Let $T$ be an ideal triangle in $(\mathcal{C},d_{\mathcal{C}})$ whose boundary $\partial \mathcal{C}$
is of regularity $\alpha$-H\"older for some $\alpha>1$. Let $a$, $b$ and $c$ be the vertices
of $T$. Let $D_a$, $D_b$ and $D_c$ be the tangent at $a$, $b$ and $c$ respectively
to $\partial \mathcal{C}$. For any two points $p$, $q$ in the plane, let $D_{pq}$ be the straight line
passing by $p$ and $q$.
Let us focus on the vertex $a$, and choose a system of coordinates in $\R^2$
such that the $x$-axes is the straight line $D_a$ and the convex $\mathcal{C}$ lies in $\R\times[0,+\infty)$.

Then Lemma \ref{lemgraph-strip} implies that for $\rho$ small enough, there is a function
$f\colon [-\rho,\rho]\to \R$ and a real number $h>0$ such that $\partial\mathcal{C}\cap([-\rho,\rho]\times[0,h])$ is the graph
of $f$. Now choose $a'\in D_{ab}$ and $a''\in D_{ac}$ such that $D_{a'a''}$ is parallel to $D_{bc}$
and $[a',a'']\subset [-\rho,\rho]\times[0,h]$.
Lemma \ref{lemfinite-volume} implies that the area of the triangle $T_a=aa'a''$ is finite.

In the same way, we built two other triangles $bb'b''$, $cc'c''$ which are of finite area.
Now the hexagon $\mathcal{H}=(a'a''b'b''c''c')$ is a compact set in $(\mathcal{C},d_\mathcal{C})$, hence of
finite area.

Thus the ideal triangle $T$ which is the union of the hexagon $\mathcal{H}$ and the triangles
$T_a$, $T_b$ and $T_c$ is of finite area.






\end{proof}

From Benoist work, mainly corollaire 1.5, a) page 184 in \cite{benoist}, 
we know that if $(\cC , d_{\cC})$ is Gromov-hyperbolic, then there is some $\alpha\in \mathopen]1,2]$
such that $\partial \mathcal{C}$ is $C^\alpha$. Hence follows  

\begin{cor}
  Let $\cC$ be a bounded open convex domain in $\R^2$ such that $(\cC , d_{\cC})$ is Gromov-hyperbolic.
Then, for any ideal triangle $T$ in $(\cC,d_{\mathcal{C}})$, we have that $\mu_{\mathcal{C}}(T)$ is finite.
\end{cor}

\subsection{Proof of Theorem \ref{thmhyp-bounded}}

The proof is done by contradiction.  

\smallskip
Assume that we can find a sequence $(\omega_n,T_n)\in T_2^\delta$ such that
$$
{\sup{\{ \mu_{\ \omega_n}(T_{n}) : n \in \N \}} = +\infty}
$$ 
and prove this is not possible.

The main idea is to use the fact that $G_2$ acts co-compactly  by isometries on the triangle-pointed
convex, to transform a converging subsequence  of $ (\omega_n,T_n)\in T_2^\delta$ into a 
sequence of convex sets $(\cC_n,\Delta)\in T_2^\delta$ evolving around a fixed ideal triangle $\Delta$.

Then in a perfect world we would be able to find a convex set $\mathcal{C}_{perfect}$
containing $\Delta$ as an ideal triangle with finite area and included in all $\mathcal{C}_n$, 
and then we would get a contradiction. 

Things are not that easy, but almost. Actually,
we will cut $\Delta$ into $4$ pieces, and for each of these pieces,
we will show that there is a convex set for which it is of finite volume
and included in $\mathcal{C}_n$ for all $n\in \N$.

\smallskip
Before going deeper into the proof, let us first make an overview of the different steps.

\begin{enumerate}[Step 1:]

\item We transform the problem to obtain a converging sequence  $(\cC_n,\Delta)\in T_2^\delta$ to $(\Omega,\Delta)$, where $\Delta$ is a fixed ideal triangle, and $\cC_n$ are  convex sets tangent to two fixed lines at two of the vertices
of $\Delta$.

\item In this step, we built a small convex set $\mathcal{G}_1\subset\cC_n$ around the vertex $(1,0)$ of
$\Delta$, which is tangent to the $x$-axe at $(1,0)$ and such that a sufficiently small section
$T_1$ of $\Delta$ containing the vertex $(1,0)$ is of finite volume $V_1$ in $\mathcal{G}_1$.

\item Reasoning as in the previous step 
we built a small convex set $\mathcal{G}_2\subset\cC_n$ around the vertex $(0,1)$ of
$\Delta$, which is tangent to the $y$-axe at $(0,1)$ and such that a sufficiently small section
$T_2$ of $\Delta$ containing the vertex $(0,1)$ is of finite volume $V_2$ in $\mathcal{G}_2$.

\item We built a small triangle $A$ which is a section of $\Delta$ admitting the vertex $(1,1)$ 
as one of its vertices and whose volume is bounded by a finite number $V_3$ for any $\cC_n$.

\item We built a convex set $\mathcal{U}$ and a compact set $S$ such that
  \begin{enumerate}
  \item for all $n$, $\mathcal{U}\subset \cC_n$.
  \item $\mu_\mathcal{U}(S)=V_4$ is finite; and
  \item $S \cup A \cup T_1 \cup T_2 =\Delta$;
  \end{enumerate}
We then conclude that for all $n$
  \begin{equation}
\begin{split}
    \mu_{\mathcal{C}_n}(\Delta)&\leq\mu_{\mathcal{C}_n}(T_1) + \mu_{\mathcal{C}_n}(T_2)  + \mu_{\mathcal{C}_n}(A) +\mu_{\mathcal{C}_n}(S) \\
&\leq  \underbrace{\mu_{\mathcal{G}_1}(T_1)}_{\leq V_1 \text{ by step 2}} +\underbrace{\mu_{\mathcal{G}_2}(T_2)}_{\leq V_2 \text{ by step 3}}
+ \underbrace{\mu_{\mathcal{C}_n}(A)}_{\leq V_3\text{ by step 4}} +\underbrace{\mu_{\mathcal{U}}(S)}_{\leq V_4 \text{ by step 5}} \\
&\leq V_1+V_2+V_3+V_4 < +\infty
\end{split}
  \end{equation}
which is absurd.
\end{enumerate}

\noindent
\istep\step\label{step0} Extraction of the subsequence.

\medskip
Following the proof of Proposition \ref{propcocompactness}, and keeping its notations, 
we can find
a sequence $(g_n)$ in $G_2$ such that $g_nT_n=\Delta$ and $g_n\omega_n=\mathcal{C}_n$,
which therefore satisfies 
\begin{equation}
\mu_{\cC_{n}}(\Delta) = \mu_{\omega_n}(T_{n})\text{.} 
\end{equation}
This implies that
\begin{equation}
\sup\{ \mu_{\cC_{n}}(\Delta) : n \in \N \} = +\infty\text{.}
\end{equation}

Furthermore, always by Proposition \ref{propcocompactness}, we have $(\mathcal{C}_n,\Delta)_{n\in \N}$ which
converges towards some $\Omega$. Recall that (for all $n\in \N$) $\mathcal{C}_n$ and $\Omega$ are
tangent to the $x$-axe at $(1,0)$, to the $y$-axe at $(0,1)$ and at $(1,1)$ to some line.

\noindent
\step\label{step1}

\medskip
We will need the following theorem

\begin{theo}\label{theouniforme}
  Let $(\cD_n)_{n\in \N}$ be a sequence of convex sets in $\R^2$ 
whose Hilbert Geometry is $\delta$-hyperbolic, for some fixed $\delta$, and
a straight line $L$. Assume that
\begin{itemize}
\item the sequence $(\cD_n)_{n\in \N}$ converges to some open convex set $\mathcal{D}$;
\item There is some $p\in L$ such that for all $n$, 
$\cD_n$ lies in the same half plane determined by $L$, and
is tangent at $p$ to $L$; 
\end{itemize}
then, taking as origin the point $p$, as $x$-axe the line $L$, and as $y$-axe an orthogonal
line to $L$, 
\begin{enumerate}[1)]
\item \label{lemuniform-func} There is a number $3a=\rho > 0$ such that for all $n \in \N$, there is a convex
function $f_{n} \colon [-3 a , 3 a] \to \R$ and  numbers  $b_{n}> 0$ and $s_n\in\R$  such that
\begin{equation}
\bD_{n} \cap \{ (x , y) \in \R^2 : x \in [-3a , 3a] \ \text{and} \ y < s_n x + b_n \} = \text{Graph}{f_{n}}.
\end{equation}
\item There is some $\mu>0$ and $\alpha>0$ (which will be made explicit in the proof) such that
$$
f_{n}(x) \leq \mu |x|^{\alpha} \quad \text{for all} \quad x \in [-a , a]\text{,}
$$
\item \label{lemlower-bound} Let 
$$ m(f_n)= \min\{f_n(-a),f_n(a)\} 
$$
then we have $u_0:=\inf \{m(f_n): n\in \N \} >0$.
\end{enumerate}
\end{theo}

We first show how to use Theorem \ref{theouniforme} to achieve the second step of the proof
of Theorem \ref{thmhyp-bounded}.
\medskip
We have just to exhibit precisely the part $T_1$ of the triangle whose area will
be bounded above, independantely of the $\delta$-hyperbolic convex set $\mathcal C_n$ we consider.

\medskip
Let 
\begin{equation}
\cD := \Omega + (-1 , 0)
\end{equation} 
(translate of $\cC$ by the vector $(-1 , 0)$) and 
\begin{equation}
\cD_{n} := \cC_{n} + (-1 , 0) \text{ for all } n \in \N\text{.}
\end{equation}
Note that since $(\Omega , d_{\Omega})$ is Gromov-hyperbolic, the same is true for $(\cD , d_{\cD})$.

\smallskip
We thus apply Theorem \ref{theouniforme} to this sequence $\cD_n$ in order to use 
Lemma \ref{lemfinite-volume} to see that a fixed triangle $T_1$ has finite area.

\smallskip
To do this, let us consider $a_{0} := (\alpha \mu)^{\frac{1}{1 - \alpha}} > 0$. 
The tangent line to $\{ (x , \mu |x|^{\alpha}) : x \in \R \}$ at the point $(-a_{0} , \mu a_{0}^{\alpha})$
is then parallel to the line $\{ (x , y) \in \R^2 : y = -x \}$.

Define 
$$
u_{1} := \min{\! \{ u_{0} , \mu a_{0}^{\alpha} \}} > 0
$$
 and pick any $u \in (0 , u_{1} / 3]$.
Applying the linear transformation of $\R^2$ given by 
$$
(x , y) \mapsto (-x (3 u / \mu)^{-1 / \alpha} , y / 3 u),
$$
we are in the situation of Lemma~\ref{lemfinite-volume} 
with
$$
\lambda= (3 u / \mu)^{1 / \alpha} / 3 u \geq 1\text{,}
$$
from which we can deduce with $\tau := 2 / 3 \in (0 , 1)$ that
the triangle
$$
\{ (x , y) \in \R^2 : x < 0\ \text{ and }\ -x < y < 2 u \}
$$
is included in the bounded open convex domain 
$$
\{ (x , y) \in \R^2 : \mu |x|^{\alpha} < 3 u \ \text{ and }\ \mu |x|^{\alpha} < y < 3 u \}
$$
and has a finite Hilbert area.

So, if we consider the triangle
$$
T_{1}(u) := \{ (x , y) \in \R^2 : x < 0 \ \text{ and }\ -x < y < 2 u \} + (1 , 0)
$$
and the bounded open convex domain
$$
\cG_{1}(u) := \{ (x , y) \in \R^2 : \mu |x|^{\alpha} < 3 u \ \text{ and }\ \mu |x|^{\alpha} < y < 3 u \} + (1 , 0),
$$
we have $T_{1}(u) \subset \cG_{1}(u)$ and $V_{1} := \mu_{\cG_{1} \! (u)}(T_{1}(u))$ is finite.

In addition, since $3 u < u_{0}$, we get that for all $n \in \N$, $\cG_{1}(u)$ is 
contained in the convex set $\cC_{n}$,
and thus $\mu_{\cC_{n}}(T_{1}(u)) \leq V_{1}$ by Proposition~\ref{propcomparison} of Appendix.

\begin{proof}[Proof of Theorem \ref{theouniforme}]
  
Let us postpone the proof of  claim 1), and prove the other two claims 

\medskip
\noindent\textbf{Claim 2)} First note that we have $f_{n} \geq 0$ and $f_{n}(0) = 0$ since $(0 , 0) \in \bD_{n}$.
In addition, as $(\cD_{n} , d_{\cD_{n}})$ is $\delta$-hyperbolic, 
Lemma~6.2, page 216, and Proposition~6.6, page 219, in Benoist \cite{benoist}
imply there is a number $H = H(\delta) \geq 1$, independant of $n$ such that
$f_{n}$ is $H$-quasi-symmetrically convex on the compact set $[-2 a , 2 a] \subset (-3 a , 3 a)$.

Therefore, by Lemma~\ref{lemHolder}, we have for $H_2=(4 H (H + 1))^{\frac{1 + a}{a}}$,
\begin{equation}
f_{n}(x) \leq 160 (H_2 + 1) M(f_{n}) |x|^{\alpha} \quad \text{for all} \quad x \in [-a , a]\text{,}
\end{equation}
where 
\begin{equation}
\alpha = 1 + \log_{2}{\!\! \Big(1 + H_2^{-1} \Big)} > 1
\end{equation}
and 
\begin{equation}
M(f_{n}) = \max\{ f_{n}(-a) , f_{n}(a) \}\text{.}
\end{equation}





We next claim that the sequence $(M(f_{n}))_{n \in \N}$ is bounded above.

Indeed, suppose that $\sup{\! \{ M(f_{n}) : n \in \N \}} = +\infty$.
If $\pi_{2} \colon \R^{2} \to \R$ denotes the projection onto the second factor, there is a number $R > 0$
such that $\pi_{2}(\cD) \subset [0 , R]$ since $\cD$ is bounded and included in $\R \times [0 , +\infty)$
(this latter point is a consequence of the fact that $\cD_{n} \subset \R \times [0 , +\infty)$ for all $n \in \N$).
Then, using $d_{H}(\pi_{2}(\bar{\cD}_{n}) , \pi_{2}(\bar{\cD})) \to 0$ as $n \to +\infty$,
there is an integer $n_{1} \in \N$ such that $\pi_{2}(\bar{\cD}_{n}) \subset [0 , 3 R]$
for all $n \geq n_{1}$.

Now there exists $n \geq n_{1}$ such that $ M(f_{n}) \geq 4 R$, that is, 
$$
\pi_{2}(-a , f_{n}(-a)) = f_{n}(-a) \geq 4 R \text{ or } \pi_{2}(a , f_{n}(a)) = f_{n}(a) \geq 4 R\text{.}
$$ 
As the points $(-a , f_{n}(-a))$ and $(a , f_{n}(a))$ are both in $\bar{\cD}_{n}$,
we get that $\pi_{2}(\bar{\cD}_{n}) \cap [4 R , +\infty) \neq \emptyset$,
which is not possible.
  
Hence there is a constant $M > 0$ such that for all $n \in \N$, we have $M(f_{n}) \leq M$,
and thus 
$$
f_{n}(x) \leq \mu |x|^{\alpha} \quad \mbox{for all} \quad x \in [-a , a],
$$
where $\mu := 160 (H + 1) M > 0$. Which proves our second claim.

\medskip
\noindent\textbf{Claim 3)} Now, for any $n \in \N$, recall that 
\begin{equation}
m(f_{n}) = \min\{ f_{n}(-a) , f_{n}(a) \}\text{.}
\end{equation}

We have to show that 
\begin{equation}
\inf \{ f_{n}(-a) : n \in \N \} \text{ and } \inf\{ f_{n}(a) : n \in \N \}
\end{equation}
are both positive numbers. So, assume that one of them,
for example the second one, is equal to zero.

Therefore there would exist a subsequence of $\Bigl(\bigl(a , f_{n}(a)\bigr)\Bigr)_{n \in \N}$ 
that converges to $(a , 0)$,
and hence $(a , 0) \in \bar{\cD}$, since $d((a , f_{n}(a)) , \cD) \to 0$ as $n \to +\infty$. 

As $(0 , 0) \in \bar{\cD}$, the whole line segment $\{ 0 \} \times [0 , a]$ would 
then be included in the convex set $\bar{\cD}$.
But $\cD$ is included in $\R \times [0 , +\infty)$ whose boundary contains $\{ 0 \} \times [0 , a]$.

This would imply that $\{ 0 \} \times [0 , a] \subset \bD$, and thus $\cD$ would not be strictly convex,
contradicting the Gromov hyperbolicity of $(\cD , d_{\cD})$ by Soci\'e-M\'ethou \cite{soth}.

\medskip
\noindent\textbf{Claim 1)} Now back to the first claim.
It suffices to use the following lemma

\begin{lema} \label{lemuniform-strip}
   There is a number $\rho > 0$ such that for all $n \in \N$, we have 
   $$
   \cD_{n} \cap ((-\infty , -\rho) \times \R) \neq \emptyset
   \quad \text{and} \quad 
   \cD_{n} \cap ((\rho , +\infty) \times \R) \neq \emptyset.
   $$
\end{lema}

Thus, given any $n \in \N$, as $\cD_{n} \subset \R \times [0 , +\infty)$,
it suffices to apply Lemma~\ref{lemgraph-strip} with $a := \rho / 3$ and $\cS = \cD_{n}$
in order to get numbers $b_{n}> 0$ and $s_n\in\R$ 
and a convex function $f_{n} \colon [-3 a , 3 a] \to \R$ such that
\begin{equation}
\bD_{n} \cap  (x , y) \in \R^2 : x \in [-3a , 3a] \ \text{and} \ y < s_n x + b_n \} = \text{Graph}{f_{n}}.
\end{equation}
\end{proof}

\begin{proof}[Proof of Lemma~\ref{lemuniform-strip}]

From the Gromov hyperbolicity of $(\cD , d_{\cD})$, we get that
the boundary $\bD$ is a $1$-dimensional submanifold of $\R^2$ of class $C^{1}$ by Karlson-Noskov 
\cite{kn}.

As $(0 , 0) \in \bD$, Lemma~\ref{lemcurve-sector} then implies that $\cD$ 
neither lies in $(0 , +\infty) \times (0 , +\infty)$,
nor in $(-\infty , 0) \times (0 , +\infty)$.

Hence, denoting by $\pi_{1} : \R^2 \to \R$ the projection onto the first factor,
$\pi_{1}(\cD)$ is an open set in $\R$ that contains $0$, 
and thus there exists a number $r > 0$ such that 
$$
[-2 r , 2 r] \subset \pi_{1}(\cD)\text{.}
$$

Since $\pi_{1}$ is continuous and 
\begin{equation}
d_{H}(\cD_{n} , \cD) \to 0 \text{ as } n \to +\infty\text{,} 
\end{equation}
we get that 
\begin{equation}
  d_{H}(\pi_{1}(\cD_{n}) , \pi_{1}(\cD)) \to 0 \text{ as } n \to +\infty\text{,}
\end{equation}
which implies there is an integer $n_{0} \in \N$ such that for all $n > n_{0}$, 
one has 
\begin{equation}
\pi_{1}(\cD_{n}) \cap (-\infty , -r) \neq \emptyset \quad \text{ and } \quad \pi_{1}(\cD_{n}) \cap (r , +\infty) \neq \emptyset\text{.}
\end{equation}

Finally, given any $n \in \{ 0 , \ldots , n_{0} \}$, 
there exists $r_{n} > 0$ such that $[-2 r_{n} , 2 r_{n}] \subset \pi_{1}(\cD_{n})$
by applying to $\cD_{n}$ the same argument as the one used for $\cD$ above.

But this implies that
\begin{equation}
  \pi_{1}(\cD_{n}) \cap (-\infty , -r_{n}) \neq \emptyset \quad \text{ and } \quad \pi_{1}(\cD_{n}) \cap (r_{n} , +\infty) \neq \emptyset\text{.}
\end{equation}

Then, choosing $\rho = \min{\{ r ,r_0,\ldots,r_{n_{0}} \}} > 0$, Lemma~\ref{lemuniform-strip} is proved.
\end{proof}

\step\label{step2}
Using the translation by the vector $(0 , -1)$ and reasoning as in Step~1 with $x$ and $y$ exchanged, 
we get numbers $\beta > 1$, $\nu > 0$, $b_{0} > 0$ and $0 < v_{1} \leq \nu b_{0}^{\beta}$ 
such that the following holds:

\begin{enumerate}
   \item The tangent line to $\{ (\nu |y|^{\beta} , y) : y \in \R \}$ at the point $(\nu b_{0}^{\beta} , b_{0})$
   is parallel to the line $\{ (x , y) \in \R^2 : y = -x \}$.
   
   \item For each $v \in (0 , v_{1} / 3]$, the triangle
   $$
   T_{2}(v) := \{ (x , y) \in \R^2 : y < 0 \text{ and } -y < x < 2 v \} + (0 , 1)
   $$
   and the bounded open convex domain
   $$
   \cG_{2}(v) := \{ (x , y) \in \R^2 : \nu |y|^{\beta} < 3 v \text{ and } \ \nu |y|^{\alpha} < x < 3 v \} + (0 , 1),
   $$
   satisfy 
$$
T_{2}(v) \subset \cG_{2}(v) \subset \cC_{n}
$$ 
for all $n \in \N$ 
and 
$$
V_{2} := \mu_{\cG_{2}(v)}\bigl(T_{2}(v)\bigr)
$$
is finite.
\end{enumerate}

\smallskip
Therefore, we deduce that $\mu_{\cC_{n}}(T_{2}(v)) \leq V_{2}$ for all $n \in \N$.

\step\label{step3} The geometric idea is similar to the two precedent steps. The only
difficulty is that the tangent line to $\mathcal C_n$ is not always the same, and we have
to be sure to make an uniform choice.

\medskip
For each $n \in \N$, consider the affine transformation $\Phi_{n}$ of $\R^2$ defined by
\begin{align*}
\Phi_{n}(1 , 1) = (0 , 0)\text{,}&& 
\Phi_{n}  (\alpha_{n} - 1 , \alpha_{n}) = (1 , 0)\text{,}&&
\Phi_{n}  (-\alpha_{n} , \alpha_{n} - 1) = (0 , 1)\text{.}
\end{align*}

As 
$$
\Phi_{n}((1 , 1) + \R_{+} (-1 , 0)) = \R_{+} (1 -\alpha_{n} , \alpha_{n})
$$
and 
$$
\Phi_{n}((1 , 1) + \R_{+} (0 , -1)) = \R_{+} (-\alpha_{n} , 1 - \alpha_{n})\text{,}
$$
since $\alpha_{n} \in [e , 1 / 2]$ (see Proposition \ref{propcocompactness}) we also have
\begin{equation} \label{equsector}
   \Phi_{n}(\Delta) \subset \{ (x , y) \in \R^2 : y \geq e |x| / (1 - e) \}\text{.}
\end{equation}

Then, applying Theorem~\ref{theouniforme} to $\Phi_{n}(\cC_{n}) \subset \R \times [0 , +\infty)$,
we get numbers $c > 0$, $\gamma > 1$ and $\kappa > 0$ and a convex function $g_{n}\colon [-c , c] \to \R$ such that 
$$
g_{n}(x) \leq \kappa |x|^{\gamma} \quad \text{for all} \quad x \in [-c , c].
$$
Next, as in claim~\ref{lemlower-bound} of Theorem~\ref{theouniforme}  in Step~\ref{step1},
 there exists a constant $w_{0} > 0$ such that
for all $n \in \N$, we have both 
\begin{align*}
g_{n}(-c) \geq w_{0} && \text{and} &&g_{n}(c) \geq w_{0}\text{.}
\end{align*}

Let
\begin{equation}
c_{0} := (e / (\kappa (1 - e)))^{1 / (\gamma - 1)} > 0\text{.}
\end{equation} 
The point $(c_{0} , \kappa c_{0}^{\gamma})$ is then the intersection point 
between the curve 
$$
\{ (x , \kappa |x|^{\gamma}) : x \in \R \}
$$
and the half line 
$$
\{ (x , y) \in \R^2 : y = e x / (1 - e), x\geq0 \}\text{.}
$$

Define
\begin{equation}
w_{1} := \min{\! \{ w_{0} , \kappa c_{0}^{\gamma} \}} > 0
\end{equation}
 and pick any $w \in (0 , w_{1} / 4]$.

Applying the linear transformation of $\R^2$ given by 
$$
(x , y) \mapsto (-x (4 w / \kappa)^{-1 / \gamma} , y / 4 w),
$$
we are in the situation of Lemma~\ref{lemfinite-volume} 
with 
$$
\lambda= \frac{e}{1-e}(4 w / \kappa)^{1 / \gamma} / 4 w.
$$

If $\lambda \geq 1$, it is an immediate application of Lemma~\ref{lemfinite-volume}. 
If $\lambda <1$, we have to do a new linear
transformation given by
$$
(x , y) \mapsto (\alpha x,\alpha^{\gamma}y)
$$
with $\alpha=\lambda^\frac{1}{1-\gamma}$, which allow to be in situation of 
Lemma~\ref{lemfinite-volume} with $\lambda = 1$

\smallskip

From this, we can deduce with $\tau := 3 / 4 \in (0 , 1)$ that
the triangle 
$$
\{ (x , y) \in \R^2 : x < 0 \text{ and } -x < y < 3 w \}
$$ 
is included in the bounded open convex domain 
$$
\{ (x , y) \in \R^2 : \kappa |x|^{\gamma} < 4 w \text{ and } \kappa |x|^{\gamma} < y < 4 w \}
$$
and has a finite Hilbert area, we denote by $V_{3}$.

So, for every $n \in \N$, if we consider the triangle
$$
A_{n}(w) := \Phi^{-1}_{n}(\{ (x , y) \in \R^2 : x < 0 \text{ and } -x < y < 3 w \})
$$
and the bounded open convex domain
$$
\cG_{n}(w) := \Phi^{-1}_{n}(\{ (x , y) \in \R^2 : \mu |x|^{\alpha} < 4 w \text{ and } \mu |x|^{\alpha} < y < 4 w \}),
$$
we have 
\begin{align*}
  A_{n}(w) \subset \cG_{n}(w)&& \text{and} &&\mu_{\cG_{n}(w)}\bigl(A_{n}(w)\bigr) = V_{3}\text{.}
\end{align*}

In addition, since $4 w < w_{0}$, we get that for all $n \in \N$, $\cG_{n}(w)$ 
is contained in the convex set $\cC_{n}$,
and thus $\mu_{\cC_{n}}(A_{n}(w)) \leq V_{3}$ by Proposition~\ref{propcomparison}.

\begin{figure}[hbtp]
  \centering
  \includegraphics[scale=.7]{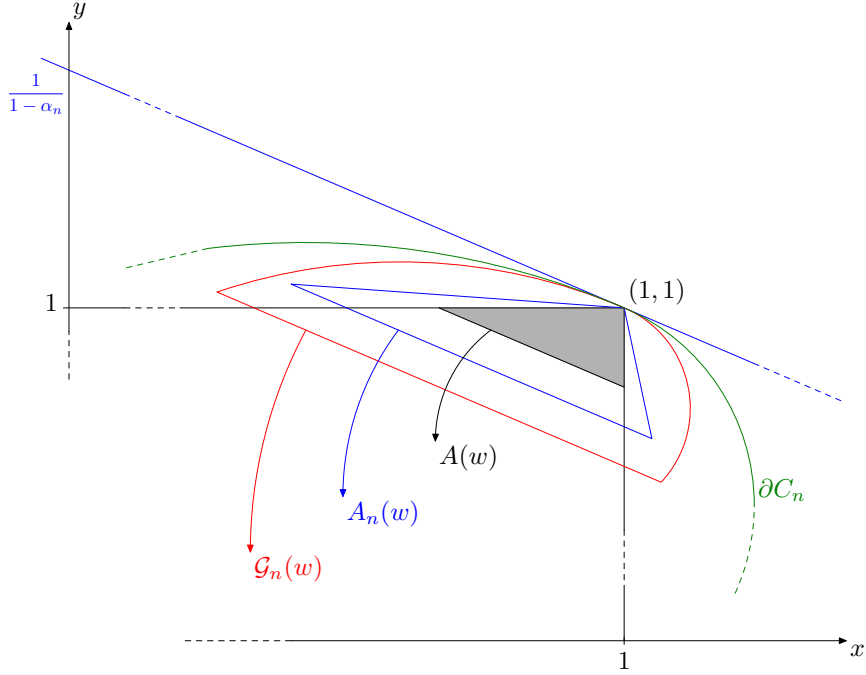}
  \caption{The triangle $A(w)$}
  \label{fig4}
\end{figure}

Now, fix $n \in \N$.

The edge of the triangle 
$$
\{ (x , y) \in \R^2 : x < 0 \mbox{ and } -x < y < 3 w \}
$$ 
that does not contain $(0 , 0)$ lies in the line $\ell := (0 , 3 w) + \R (1 , 0)$.
Hence, the edge of the triangle $A_{n}(w)$ that does not contain $(1 , 1)$ 
lies in the line 
\begin{multline}
\ell_{n} := \Phi^{-1}_{n}(\ell) 
= 
\Phi^{-1}_{n}(0 , 3 w) + \R ( \Phi_{n})^{-1}(1 , 0)\\
= (1 - 3 \alpha_{n} w , 1 + 3 (\alpha_{n} - 1) w) + \R (\alpha_{n} - 1 , \alpha_{n}).
\end{multline}
The $x$-coordinate $x_{n}$ of the intersection point between $\ell_{n}$ and the line 
$$
\{ (x , y) \in \R^2 : y = 1 \}
$$
is then equal to 
$$
x_{n} = 1 - 3 \alpha_{n} w + s (\alpha_{n} - 1)
$$ 
with $s = 3 (1 - \alpha_{n}) w / \alpha_{n}$.
From $\alpha_{n} \in [e , 1 / 2]$, we get that $s > 0$, and thus 
\begin{equation}
x_{n} < 1 - 3 \alpha_{n} w < 1 - 3 e w\text{.}
\end{equation}

On the other hand, the $y$-coordinate $y_{n}$ of the intersection point between $\ell_{n}$ 
and the line 
$$
\{ (x , y) \in \R^2 : x = 1 \}
$$
is equal to 
$$
y_{n} = 1 + 3 (\alpha_{n} - 1) w + t \alpha_{n}
$$ with $t = 3 \alpha_{n} w / (\alpha_{n} - 1)$.
As $\alpha_{n} \in [e , 1 / 2]$, we have $t < 0$, and thus 
\begin{equation}
y_{n} < 1 + 3 (\alpha_{n} - 1) w < 1 - w < 1 - 3 e w\text{.}
\end{equation}

Using Equation~\ref{equsector}, this proves that the fixed triangle $A(w)$
whose vertices are the points $(1 , 1)$, $(1 - 3 e w , 1)$ and $(1 , 1 - 3 e w)$
is included in the triangle $A_{n}(w)$.

Conclusion: by Proposition~\ref{propcomparison}, $\mu_{\cC_{n}}(A(w)) \leq V_{3}$ for all $n \in \N$.

\begin{figure}[hbtp]
  \centering
  \includegraphics[scale=.7]{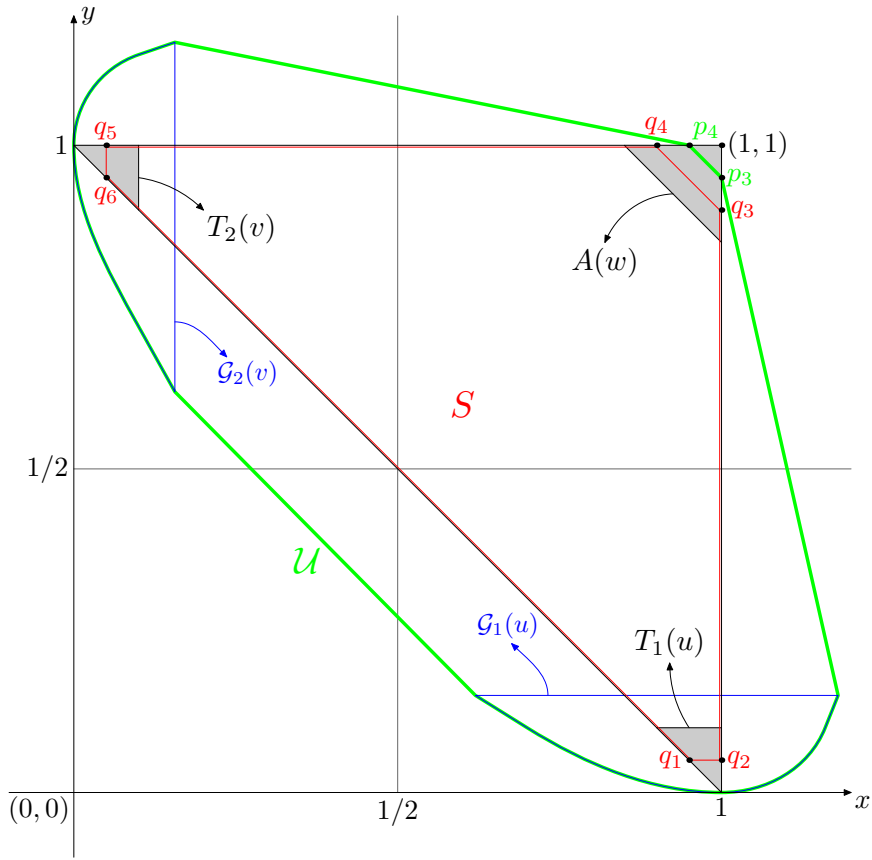}
  \caption{The convex set $\mathcal{U}$ and the compact set $S$}
  \label{fig5}
\end{figure}

\step\label{step4}
Let us introduce the points 
\begin{align*}
p_{3} &:= (1 , 1 - e w)\text{,}&
p_{4} &:= (1 - e w , 1)\text{,}&
\end{align*}
by step \ref{step3},  $p_3$ and $p_4$ are in $A(w)$ (see figure \ref{fig5}).

Now consider $\mathcal{U}$ the convex hull of $\cG_{1}(u)\cup\cG_2(v)\cup\{p_3,p_4\}$,
recall that for all $n \in \N$ we have
\begin{eqnarray*}
 \cG_{1}(u) &\subset& \cC_{n}\text{, by step \ref{step1};}\\
\cG_{2}(v) &\subset& \cC_{n}\text{, by step \ref{step2};}\\
p_{3} , p_{4} \in A(w) &\subset& \cC_{n}\text{, by step \ref{step3},}
\end{eqnarray*}
hence the fixed bounded open convex domain $\mathcal{U}$ is included in 
the convex set $\cC_{n}$, for all $n\in \N$.

Now, if $S$ is the closed convex hull in $\R^2$ of the points
\begin{align*}
  q_{1} &:= (1 - u , u)\text{,}&
q_{2} &:= (1 , u)\text{,}&
q_{3} &:= (1 , 1 - 2 e w)\text{,}\\
q_{4} &:= (1 - 2 e w , 1)\text{,}&
q_{5} &:= (v , 1)\text{,}&
q_{6} &:= (v , 1 - v)\text{,}
\end{align*}
we have $S \in \mathcal{U}$, and thus $V_{4} := \mu_{\mathcal{U}}(S)$ is finite, since $S$ is compact.

This gives that $S \subset \cC_{n}$ with $\mu_{\cC_{n}}(S) \leq V_{4}$ for all $n \in \N$.

Finally, as $\Delta \subset T_{1}(u) \cup T_{2}(v) \cup A(w) \cup S$, we get 
\begin{multline}
\mu_{\cC_{n}}(\Delta) 
\leq 
\mu_{\cC_{n}}(T_{1}(u)) + \mu_{\cC_{n}}(T_{2}(v))\\ + \mu_{\cC_{n}}(A(w)) + \mu_{\cC_{n}}(S)
\\
\leq
V_{1} + V_{2} + V_{3} + V_{4} =: V < +\infty\text{.}  
\end{multline}

But this is in contradiction with the assumption 
$$
\sup{\{ \mu_{\cC_{n}}(\Delta) : n \in \N \}} = +\infty\text{,}
$$
hence Theorem~\ref{thmhyp-bounded} is proved.
\appendix

\section{Technical lemmas}

We recall, without proof, Proposition~5 in \cite{cvv} page 208.

\begin{prop} \label{propcomparison}
  Let $({\mathcal A},d_{\mathcal A})$ and $({\mathcal B},
  d_{\mathcal B})$ be two Hilbert's geometries such that ${\mathcal A} \subset {\mathcal B} \subset \mathbb{R}^{n}$. Then~:
    
    \begin{enumebis}
    \item \label{v1} The Finsler metrics $F_{\mathcal A}$ and $F_{\mathcal B}$ 
      satisfy $F_{\mathcal B}(\,p,v) \leq F_{\mathcal A}(\,p,v)$
      for all $p\in \mathcal A$ and all none null $v \in \mathbb R^{n}$
      with equality, if, and only if $p_\mathcal{A}^{-} =
      p_\mathcal{B}^{-}$ and $p_\mathcal{A}^{+} = p_\mathcal{B}^{+}$.
       
    \item \label{v2}  If $p, q \in \mathcal{A}$, we have
      $d_\mathcal{B}(\,p,q) \leq d_\mathcal{A}(\,p,q)$.
       
    \item \label{v3} For all $p \in \mathcal{A}$, we have $\mu \left(
        B_\mathcal{A}(\,p) \right) \leq \mu \left( B_\mathcal{B}(\,p) \right)$ 
        with equality, if, and only if $\mathcal{A} = \mathcal{B}$.
       
    \item \label{v4} For any Borel set $A$ in $\mathcal{A}$, we have
      $\mu_{\mathcal{B}}(A) \leq \mu_{\mathcal{A}}(A)$ with equality, if, and only if $\mathcal{A} = \mathcal{B}$.
   \end{enumebis}
\end{prop}


\begin{lema} \label{lemcurve-sector}
   Fix $s \in \R$ and consider
   the half closed cone $C = \{ (x , y) \in \R^{2} : x \geq 0 \mbox{ and } y \leq s x \}$ in $\R^{2}$.
   Then we have:
   
   \begin{enumerate}[1)]
      \item For any $\varepsilon > 0$ and any parameterized curve $\sigma : (-\varepsilon, \varepsilon) \to \R^{2}$
      that is differentiable at $t = 0$, if $\sigma(0) = (0 , 0)$ and $\sigma((-\varepsilon , \varepsilon)) \subset C$,
      then $\ \sigma'(0) = (0 , 0)$.

      \item 
      For any $1$-dimensional topological submanifold $\Gamma$ of $\R^{2}$,
      if $(0 , 0) \in \Gamma$ and $\Gamma \subset C$, 
      then $\Gamma$ is not a differentiable submanifold of $\R^2$ at $(0 , 0)$.
   \end{enumerate}
\end{lema}

\begin{proof}
\textbf{Point 1:}
Let $\ \sigma(t) = (x(t) , y(t))$ for all $t \in (-\varepsilon,\varepsilon)$.
As $x(t) \geq 0 = x(0)$ for all $t \in (-\varepsilon , \varepsilon)$, the function $x \colon  (-\varepsilon , \varepsilon) \to \R$ 
has a local minimum at $t = 0$, and thus $x'(0) = 0$.

On the other hand, for all $t \in (-\varepsilon , \varepsilon)$, we have $y(t) \leq s x(t)$,
or equivalently $y(t) - s x(t) \leq 0 = y(0) - s x(0)$.
This shows that the function $y - s x\colon  (-\varepsilon , \varepsilon) \to \R$
has a local minimum at $t = 0$, and thus $(y - s x)'(0) = 0$.
But $x'(0) = 0$ and hence $y'(0) = 0$, which proves the first point
of the lemma.

\medskip
\noindent
\textbf{Point 2:}
Assume that $\Gamma$ is a $1$-dimensional differentiable submanifold of $\R^{2}$ at $(0 , 0)$.

Then we can find open sets $U$ and $V$ in $\R^2$ that contain $(0 , 0)$
together with a diffeomorphism $\Phi \colon  U \to V$ satisfying 
$\Phi(U \cap \Gamma) = V \cap (\R \times \{ 0 \})$ and $\Phi(0 , 0) = (0 , 0)$.

Let $\varepsilon > 0$ such that $(-\varepsilon , \varepsilon) \times \{ 0 \} \subset V \cap (\R \times \{ 0 \})$,
and consider the parameterized curve $\sigma : (-\varepsilon , \varepsilon) \to \R^{2}$ defined by
$\sigma(t) = \Phi^{-1}(t , 0)$.

As $\sigma$ is differentiable at $t = 0$ and satisfies 
$\ \sigma((-\varepsilon , \varepsilon)) \subset U \cap \Gamma\subset \Gamma \subset C$ and $\ \sigma(0) = (0 , 0)$,
we get from Point~1) above that $\ \sigma'(0) = (0 , 0)$, which implies that
$(\Phi \circ \sigma)'(0) = (0 , 0)$ by the chain rule.

But a direct calculation gives $(\Phi \circ \sigma)(t) = (t , 0)$ for all $t \in (-\varepsilon , \varepsilon)$,
and hence $(\Phi \circ \sigma)'(0) = (1 , 0) \neq (0 , 0)$.

So $\Gamma$ cannot be a $1$-dimensional differentiable submanifold of $\R^{2}$ at $(0 , 0)$,
proving the second point of the lemma.

\end{proof}
 
\begin{lema}  \label{lemgraph-strip}
   Let $\rho > 0$ and $\cS$ be an open convex domain in $\R^2$ 
   that lies in $\R \times (0 , +\infty)$. \\
   If $\cS \cap ((-\infty , -\rho) \times \R) \neq \emptyset$ 
   and
   $\cS \cap ((\rho , +\infty) \times \R) \neq \emptyset$,
   then there exist $s \in \R$, $b > 0$ and a function $f \colon [-\rho , \rho] \to \R$ such that
   $$
   \bS \cap \{ (x , y) \in \R^2 : x \in [-\rho , \rho] \ \text{and} \ y < s x + b \} = {\rm Graph}{f}\text{.}
   $$
\end{lema}

\begin{proof}
Pick 
$$
p_{0} = (x_{0} , y_{0}) \in \cS \cap ((-\infty , -\rho) \times \R)
$$
and 
$$
p_{1} = (x_{1} , y_{1}) \in \cS \cap ((\rho , +\infty) \times \R)\text{.}
$$
The closed line segment $L$ with vertices $p_{0}$ and $p_{1}$ 
then lies in the convex set $\cS$.
Denoting by $\pi_{1} \colon  \R^2 \to \R$ the projection onto the first factor, 
we thus get $[-\rho , \rho] \subset [x_{0} , x_{1}] = \pi_{1}(L) \subset \pi_{1}(\cS)$. 
Since $\cS \subset \R \times (0 , +\infty)$, this allows us to consider
the function $f \colon [-\rho , \rho] \to \R$ defined by
$$
f(x) = \inf{\{ y \geq 0 : (x , y) \in \bar{\cS} \}}. 
$$
Fix $x \in [-\rho , \rho]$.

Given any $z \geq 0$ such that $(x , z) \in \bar{\cS}$, we have by compactness
$$
f(x) = \inf{\{ y \in [0 , z] : (x , y) \in \bar{\cS} \}}
= \min{\{ y \in [0 , z] : (x , y) \in \bar{\cS} \}}\text{,}
$$
and thus $(x , f(x)) \in \bar{\cS}$.

If $(x , f(x))$ were in $\cS$, there would exist
$\varepsilon > 0$ such that 
$$
[x - \varepsilon , x + \varepsilon] \times [f(x) - \varepsilon , f(x) + \varepsilon] \subset \cS \subset \R \times [0 , +\infty)\text{,}
$$
and thus we would get $f(x) - \varepsilon \in {\{ y \geq 0 : (x , y) \in \bar{\cS} \}}$.
But this contradicts the very definition of $f(x)$.
Therefore, we have $(x , f(x)) \in \bS$.

\medskip

Now let $s = (y_{1} - y_{0}) / (x_{1} - x_{0})$ and $b = (x_{1} y_{0} - x_{0} y_{1}) / (x_{1} - x_{0}) > 0$.
The equation of the straight line containing $L$ is then $y = s x + b$.

\medskip

Since for all $x \in [-\rho , \rho]$, the point $(x , s x + b) \in L \subset \bar{\cS}$, we get 
$f(x) \leq s x + b$ from the definition of $f$. As $(x , f(x)) \in \bS$ and $L \cap \bS = \emptyset$, 
we also have $f(x) \neq s x + b$.
Hence 
$$
\text{Graph}{f} \subset \bS \cap \{ (x , y) \in \R^{2} : x \in [-\rho , \rho] \ \text{and} \ y < s x + b \}\text{.}
$$

\medskip

On the other hand, for any given 
$(x , z) \in \bS \cap \{ (x , y) \in \R^{2} : x \in [-\rho , \rho] \ \text{and} \ y < s x + b \}$, 
assume there is $y \geq 0$ with $(x , y) \in \bar{\cS}$ satisfying $y < z$.
Then $(x , z)$ is in the triangle whose vertices are $p_{0}$, $p_{1}$ and $(x , y)$, 
which is not possible since this triangle lies in $\cS$ 
(the interior of the closure of a convex set in $\R^{n}$ is equal 
to the interior of that convex set in $\R^{n}$)
and $(x , z) \in \bS$. Therefore $z \leq y$, which shows that $z = f(x)$ by the definition of $f$.
This proves that 
$$
\bS \cap \{ (x , y) \in \R^{2} : x \in [-\rho , \rho] \ \text{and} \ y < s x + b \} \subset \text{Graph}{f}\text{.}
$$
\end{proof}

\begin{rmq} \label{remfconvex}
 The function $f$ obtained in Lemma~\ref{lemgraph-strip} satisfies $f \geq 0$ and is automatically convex 
since its epigraph is equal to the convex set in $\R^{2}$ equal to the union of the convex set 
$$
\bar{\cS} \cap \{ (x , y) \in \R^{2} : x \in [-\rho , \rho] \ \text{and} \ y < s x + b \} \subset \R^{2}
$$
(intersection of two convex sets in $\R^{2}$) and the convex set 
$$
\{ (x , y) \in \R^{2} : x \in [-\rho , \rho] \ \text{and} \ y \geq s x + b \} \subset \R^{2}\text{.}
$$
\end{rmq}

\begin{lema} \label{lemfinite-volume}
   Let $\alpha > 1$, $\lambda\geq 1$ and $\tau \in (0 , 1)$. 
   Consider the bounded open convex domain 
   \begin{equation}
     \cG = \{ (x , y) \in \R^2 : -1 < x < 1 \mbox{ and } |x|^{\alpha} < y < 1 \}     
   \end{equation}
   and the triangle 
\begin{equation}
T= \{ (x , y) \in \R^2 : x > 0 \mbox{ and } \lambda x < y < \tau \}\text{.}
\end{equation}
   Then we have $T \subset \cG$ and the area $\mu_{\cG}(T)$ is finite.
\end{lema}

\begin{figure}[hbtp]
  \centering
  \includegraphics{dessins.1}
  \caption{figure 1}
  \label{fig1}
\end{figure}

\begin{proof}
\istep
 \step
For each $p = (x , y) \in T$, let $B_{\cG}(p) = \{ v \in \R^{2} : F_{\cG}(p , v) < 1 \}$
be the open unit ball in $T_{p}\cG = \R^{2}$ of the norm $F_{\cG}(p , \cdot)$.

An easy computation shows that the vectors 
$$
v_{1} = ((y^{2 / \alpha} - x^{2}) / y^{1 / \alpha} , 0)
\quad \text{and} \quad
v_{2} = (0 , 2 (1 - y) (y - x^{\alpha}) / (1 - x^{\alpha}))
$$
are in the boundary $\partial B_{\cG}(p)$ of $B_{\cG}(p)$.

As $B_{\cG}(p)$ is convex and symmetric about $(0 , 0)$ in $T_{p}\cG = \R^2$,
we get that the rhombus defined as the convex hull of $v_{1}$, $v_{2}$, $-v_{1}$ and $-v_{2}$ 
is included in the closure of $B_{\cG}(p)$ in $T_{p}\cG = \R^2$.

Therefore the euclidien volume of this rhombus is less than or equal to that of $B_{\cG}(p)$,
which writes 
$$
\vol\bigl(B_{\cG}(p)\bigr) \geq 4 \frac{(1 - y) (y - x^{\alpha}) (y^{2 / \alpha} - x^{2})}{y^{1 / \alpha} (1 - x^{\alpha})}~.
$$
Since $1 - x^{\alpha} \leq 1$ and $1 - y \geq 1 - \tau$, we then deduce
$$
\vol\bigl(B_{\cG}(p)\bigr) \geq 4 (1 - \tau) \frac{(y - x^{\alpha}) (y^{2 / \alpha} - x^{2})}{y^{1 / \alpha}}~.
$$

\medskip

\step
From the inequality obtained in Step~1, we have
\begin{equation} \label{equ:volume}
   \mu_{\cG}(T) = \pi \!\! \int \!\!\! \int_{T} \frac{\mathrm{d} x\, \mathrm{d} y}{\vol{B_{\cG}(p)}} 
   \leq \frac{\pi}{4 (1 - \tau)} I,
\end{equation}
where 
$$
\disp I := \int \!\!\! \int_{T} \frac{y^{1 / \alpha} \mathrm{d} x \,\mathrm{d} y}{(y - x^{\alpha}) (y^{2 / \alpha} - x^{2})}
\text{.}
$$

Now, using the change of variables $\Phi \colon  (0 , +\infty)^{2} \to (0 , +\infty)^{2}$
defined by 
$$
(s , t) = \Phi(x , y) := (x / y^{1 / \alpha} , x)\text{,}
$$
whose Jacobian at any ${(x , y) \in (0 , +\infty)^{2}}$ is equal to $x/(\alpha y^{1 + 1 / \alpha})$,
we get
$$
I = \alpha \!\! \int \int_{\Phi(T)} \frac{\mathrm{d} s\, \mathrm{d} t}{t (1 - s^{\alpha}) (1 - s^{2})}
$$ 
with ${\Phi(T)=\{(s ,t) \in \R^2 :0 < t < \tau / \lambda\text{ and } t\cdot\tau^{-1 / \alpha} < s < \lambda^{-1 / \alpha}\cdot t^{1 - 1 / \alpha} \}\text{.}}$

So,
%
\begin{eqnarray*}
   I & = & \alpha \int_{0}^{\tau / \lambda}{\frac{1}{t} 
   \left( \int_{\tau^{-1 / \alpha} t}^{\lambda^{-1 / \alpha} t^{1 - 1 / \alpha}}
   {\frac{1}{(1 - s^{\alpha}) (1 - s^{2})}}{\mathrm{d}s} \right)}{\mathrm{d}t} \\
   & \leq & \alpha \int_{0}^{\tau / \lambda}{\frac{1}{t} 
   \left( \int_{\tau^{-1 / \alpha} t}^{\lambda^{-1 / \alpha} t^{1 - 1 / \alpha}}
   {\frac{1}{(1 - \tau^{\alpha - 1} \lambda^{-\alpha}) (1 - \tau^{2 - 2 / \alpha} \lambda^{-2})}}{\mathrm{d}s} \right)}{\mathrm{d}t},
\end{eqnarray*}
since $(s , t) \in \Phi(T)$ implies
$$
1 - s^{\alpha} \geq 1 - t^{\alpha - 1} / \lambda\geq \tau^{\alpha - 1} \lambda^{-\alpha} 
\quad \mbox{and} \quad 
1 - s^{2} \geq 1 - \lambda^{-2 / \alpha} t^{2 - 2 / \alpha} \geq \tau^{2 - 2 / \alpha} \lambda^{-2}.
$$
Therefore, one has
\begin{eqnarray*}
   I & \leq & \alpha \int_{0}^{\tau / \lambda}{\frac{\lambda^{-1 / \alpha} t^{-1 / \alpha} - \tau^{-1 / \alpha}}
   {(1 - \tau^{\alpha - 1} \lambda^{-\alpha}) (1 - \tau^{2 - 2 / \alpha} \lambda^{-2})}}{\mathrm{d}t} \\
   & \leq & \Lambda\int_{0}^{\tau / \lambda}{\frac{1}{t^{1 / \alpha}}}{\mathrm{d}t},
\end{eqnarray*}
where $\displaystyle{ \Lambda:= \frac{\alpha \lambda^{-1 / \alpha}}{(1 - \tau^{\alpha - 1} \lambda^{-\alpha}) (1 - \tau^{2 - 2 / \alpha} \lambda^{-2})}}$~.

Since $1 / \alpha < 1$, this shows that $I < +\infty$, and Equation~\ref{equ:volume} proves the lemma.

\end{proof}

\begin{defi}
   Given a number $K \geq 1$ and an interval $I \subset \R$, a function $f\colon I \to \R$
   is said to be $K$-quasi-symmetric if and only if one has:
\begin{multline}
   \forall x \in I,~\forall h \in \R,~(x + h \in I \text{ and } x - h \in I) \\
   \imp |f(x + h) - f(x)| \leq K |f(x) - f(x - h)|\text{.} 
\end{multline}
\end{defi}

\begin{defi}
   Given a number $H \geq 1$ and an interval $I \subset \R$, a function $f\colon I \to \R$
   is said to be $H$-quasi-symmetrically convex if and only if 
   it is convex, differentiable and has the following property:
\begin{multline}
   \forall x \in I,~\forall h \in \R,~(x + h \in I \mbox{ and } x - h \in I)\\ \imp D_{x}(h) \leq H D_{x}(-h), 
\end{multline}
   where 
$$
D_{x}(h) := f(x + h) - f(x) - f'(x) h\text{.}
$$
\end{defi}

\begin{lema} \label{lemHolder}
   Let $a > 0$, $H \geq 1$ and $f\colon [-2 a , 2 a] \to \R$ a $H$-quasi-symmetrically convex function
   that satisfies $f \geq 0$ and $f(0) = 0$.
   Define
   \begin{equation}
    H_2= \bigl(4 H (H + 1)\bigr)^{\frac{1 + a}{a}} >1
   \end{equation}
and
\begin{equation}
\alpha = 1 + \log_{2}\, ( 1+ H_2^{-1}) > 1  
\end{equation}
   and  $M(f) = \max\{ f(-a) , f(a) \}$.
   Then we have
\begin{equation}
   f(x) \leq 160 (H_2 + 1) M(f) |x|^{\alpha} \quad \text{for all} \quad x \in [-a , a].
\end{equation}
\end{lema}

Before proving this lemma, recall the two following results due to Benoist \cite{benoist}.

\bigskip

\begin{lema}[\cite{benoist}, Lemma~5.3.b), page 204] \label{lemBenoist1}
   Let $a > 0$, $H \geq 1$ and $f \colon [-2 a , 2 a] \to \R$ a $H$-quasi-symmetrically convex function.
   Then the restriction of the derivative $f'$ to $[-a , a]$ is $K$-quasi-symmetric,
   where $K = (4 H (H + 1))^{\frac{1 + a}{a}} \geq 1$.
\end{lema}

\bigskip

\begin{lema}[\cite{benoist}, Lemma~4.9.a), page 203] \label{lemBenoist2}
   Let $a > 0$, $K \geq 1$ and $f \colon [-a , a] \to \R$ a differentiable convex function. \\
   If the derivative $f'$ is $K$-quasi-symmetric,
   then for all $x , y \in [-a , a]$, we have
   $$
   |f'(x) - f'(y)| \leq 160 (1 + K) ||f||_\infty |x - y|^{\alpha - 1},
   $$
   where $\alpha = 1 + \log_{2}{\! \big( 1 + 1 / K \big)} > 1$.
\end{lema}

\begin{proof}[Proof of lemma \ref{lemHolder}]
Using Lemma~\ref{lemBenoist1}, we have that the derivative $f'$
of $f$ is $K$-quasi-symmetric when restricted to $[-a , a]$.
But then, according to Lemma~\ref{lemBenoist2} with $y = 0$ and the fact that $f'(0) = 0$
since $0$ is the minimum of $f$, we get that
$$
|f'(x)| \leq 160 (1 + K) \max_{t \in [-a , a]}{\! |f(t)|} \, |x|^{\alpha - 1}
\quad \mbox{for all} \quad x \in [-a , a].
$$
Now, the convexity of $f$ implies that $f'$ is a non-decreasing function.
As $f'(0) = 0$, we have that $f'(x) \leq 0$ for all $x \in [-a, 0]$
and $f'(x) \geq 0$ for all $x \in [0 , a]$.
Hence, $f$ is a function that is non-increasing on $[-a , 0]$ 
and non-decreasing on $[0 , a]$, which yields to
$\max_{t \in [-a , a]}{\! |f(t)|} = M(f)$ and
\begin{equation} \label{equderivative}
   |f'(x)| \leq 160 (1 + K) M(f) |x|^{\alpha - 1}
   \quad \mbox{for all} \quad x \in [-a , a].
\end{equation}
Choosing an arbitrary $u \in [-a , a]$ and applying Taylor's theorem to $f$
between $0$ and $u$, we get the existence of $\vartheta \in (0 , 1)$
such that 
$$
f'(\vartheta u) u = f(u) - f(0) = f(u).
$$
Therefore, plugging $x = \vartheta u \in [-a , a]$ in Equation~\ref{equderivative}
and multiplying by $|u|$, one has
$$
|f(u)| = |f'(\vartheta u)| |u| \leq 160 (1 + K) M(f) |\vartheta u|^{\alpha - 1} |u| 
\leq 160 (1 + K) M(f) |u|^{\alpha}
$$
since $|\vartheta u| \leq |u|$. This proves Lemma~\ref{lemHolder}.
\end{proof}

\nocite{*}
\bibliographystyle{amsalpha}
\bibliography{airhyp}

\providecommand{\bysame}{\leavevmode\hbox to3em{\hrulefill}\thinspace}
\providecommand{\MR}{\relax\ifhmode\unskip\space\fi MR }
\providecommand{\MRhref}[2]{%
  \href{http://www.ams.org/mathscinet-getitem?mr=#1}{#2}
}
\providecommand{\href}[2]{#2}
\begin{thebibliography}{CVV04}

\bibitem[BBI01]{bbi}
D.~Burago, Y.~Burago, and S.~Ivanov, \emph{A course in metric geometry},
  Graduate Studies in Mathematics, vol.~33, American Mathematical Society,
  2001.

\bibitem[Ben60]{benzecri}
J.-. Benz{\'e}cri, \emph{Sur les vari\'et\'es localement affines et localement
  projectives}, Bull. Soc. Math. France \textbf{88} (1960), 229--332.
  \MR{MR0124005 (23 \#A1325)}

\bibitem[Ben03]{benoist}
Y.~Benoist, \emph{Convexes hyperboliques et fonctions quasisym\'etriques},
  Publ. Math. Inst. Hautes \'Etudes Sci. (2003), no.~97, 181--237. \MR{2 010
  741}

\bibitem[BH99]{bridsonh}
M.~R. Bridson and A.~Haefliger, \emph{Metric spaces of non-positive curvature},
  Grundlehren der mathematischen Wissenshaften, vol. 319, Springer, 1999, A
  Series of Comprehensive Studies in Mathematics.

\bibitem[CV]{cv1}
B.~Colbois and C.~Vernicos, \emph{Bas du spectre et delta-hyperbolicit{\'e} en
  g{\'e}om{\'e}trie de {H}ilbert plane}, to appear in Bulletin de la SMF.

\bibitem[CV04]{cv}
B.~Colbois and P.~Verovic, \emph{Hilbert geometry for strictly convex domains},
  Geom. Dedicata \textbf{105} (2004), 29--42. \MR{MR2057242 (2005e:53111)}

\bibitem[CVV04]{cvv}
B.~Colbois, C.~Vernicos, and P.~Verovic, \emph{L'aire des triangles id{\'e}aux
  en g{\'e}om{\'e}trie de {H}ilbert}, L'enseignement math{\'e}matique
  \textbf{50} (2004), no.~3-4, 203--237.

\bibitem[dlH93]{dlharpe}
P.~de~la Harpe, \emph{On {H}ilbert's metric for simplices}, Geometric group
  theory, Vol.\ 1 (Sussex, 1991), Cambridge Univ. Press, Cambridge, 1993,
  pp.~97--119.

\bibitem[Hil71]{dhilbert}
D.~Hilbert, \emph{Les fondements de la g{\'e}om{\'e}trie, edition critique
  pr{\'e}par{\'e} par p.~rossier}, Dunod, 1971, (voir Appendice I).

\bibitem[KN02]{kn}
A.~Karlsson and G.~A. Noskov, \emph{The {H}ilbert metric and {G}romov
  hyperbolicity}, Enseign. Math. (2) \textbf{48} (2002), no.~1-2, 73--89.
  \MR{2003f:53061}

\bibitem[Sam88]{ps}
P.~Samuel, \emph{Projective geometry}, Undergraduate texts in Mathematics,
  Springer, 1988, also available in french.

\bibitem[SM00]{soth}
{\'E}.~Soci{\'e}-M{\'e}thou, \emph{Comportement asymptotiques et rigidit{\'e}s
  en g{\'e}om{\'e}tries de hilbert}, Th{\`e}se de doctorat, Universit{\'e} de
  Strasbourg, 2000,
  http://www-irma.u-strasbg.fr/irma/publications/2000/00044.ps.gz.

\bibitem[SM02]{so}
\bysame, \emph{Caract\'erisation des ellipso\"\i des par leurs groupes
  d'automorphismes}, Ann. Sci. \'Ecole Norm. Sup. (4) \textbf{35} (2002),
  no.~4, 537--548. \MR{1 981 171}

\bibitem[SM04]{so2}
\bysame, \emph{Behaviour of distance functions in {H}ilbert-{F}insler
  geometry}, Differential Geom. Appl. \textbf{20} (2004), no.~1, 1--10.
  \MR{2004i:53112}

\end{thebibliography}

\end{document}